# Teichmüller geodesics that do not have a limit in $\mathcal{PMF}$

ANNA LENZHEN

We construct a Teichmüller geodesic which does not have a limit on the Thurston boundary of the Teichmüller space. We also show that for this construction the limit set is contained in a one-dimensional simplex in $\mathcal{PMF}$.

## 1 Introduction

In this paper we consider a problem in Teichmüller geometry at infinity. Recall that Teichmüller space $\mathcal{T}_g$ of a closed oriented surface of genus $g$ equipped with Teichmüller metric is a complete geodesic metric space. There is a natural way to compactify $\mathcal{T}_g$, by fixing a base point and considering all geodesic rays through that point. However Kerckoff proved that the Teichmüller, or visual, compactification was base point dependent by showing that the action of the mapping class group on $\mathcal{T}_g$ does not extend to the boundary. Another natural compactification is due to Thurston, by projective measured foliations $\mathcal{PMF}$, to which the action of the mapping class group does extend.

Teichmüller geodesic rays are associated with quadratic differentials which in turn are closely related to measured foliations: geodesics are described by scaling horizontal and vertical measured foliations of a quadratic differential. Therefore it is natural to compare the two compactifications, in particular, to study the behavior of a Teichmüller ray with respect to the Thurston compactification. The question is quite nontrivial. A geodesic is defined by deforming the flat metric. On the other hand, the Thurston compactification is defined using the hyperbolic metric, and there is no obvious way to compare the two metrics. Masur [8] proved that in almost every direction through every point, geodesic rays have a limit on the boundary. In particular Masur considered geodesic rays given by quadratic differentials with uniquely ergodic vertical measured foliation and showed that a ray defined by a differential $q$ with vertical measured foliation $F$ converges to the class of $F$. Quadratic differentials with closed vertical trajectories were considered in the same paper, and it was shown that corresponding geodesic rays also converge. However the limit in this case might not be the class of $F$ itself. It turned out that the ray approaches the barycenter of the simplex of measures on $F$, ie in the limit all the leaves have the same weight.

In this paper we give a first example of a geodesic ray which does not have a limit on the Thurston's boundary. Our main result is the following:

**Theorem 1** *There exists a Teichmüller geodesic ray which does not converge in Thurston's compactification.*

The proof is by construction of a divergent ray, which is found by gluing two copies of a square torus cut along slits whose slopes satisfy certain conditions, to obtain a surface of genus 2, and then scaling the vertical and horizontal foliations of the corresponding quadratic differential. We also describe the limit set of the divergent ray. More precisely, we show that it is contained in a one dimensional simplex on the boundary.

## 2 Preliminaries

We refer the reader to Imayoshi and Tanaguchi [4] and *Travaux* [3] for more information on Teichmüller theory. Let $M$ be a closed surface of genus $g \geq 2$. Recall that the *Teichmüller space* $\mathcal{T}_g$ is the space of equivalence classes of conformal structures $X$ on $M$. The equivalence relation is defined by considering two structures $X_1$ and $X_2$ equivalent if there is a biholomorphic map from $X_1$ to $X_2$ which is isotopic to the identity on $M$.

Let $X_1$ and $X_2$ be two points in $\mathcal{T}_g$. The *Teichmüller distance* between $X_1$ and $X_2$ is defined to be
$$d(X_1, X_2) = \frac{1}{2} \log K$$
where $K$ is the smallest number such that there is a homeomorphism homotopic to the identity on $M$ which is a $K$–quasiconformal map between $X_1$ and $X_2$. There is a unique quasiconformal map from $X_1$ to $X_2$ realizing this distance, called *Teichmüller mapping*.

A *holomorphic quadratic differential* $q$ (see Strebel [11] for details) on a Riemann surface $X$ is an assignment to each chart $(U_\alpha, z_\alpha)$ of $X$ a holomorphic function $q_{z_\alpha}$ with the property
$$q_{z_\beta}(z_\beta) \left(\frac{dz_\alpha}{dz_\beta}\right)^2 = q_{z_\alpha}(z_\alpha)$$
in $U_{z_\alpha} \cap U_{z_\beta}$. The norm or area of $q$ is defined by $\|q\| = \int_X |q(z)||dz|^2$. The vector space $Q_0$ of all holomorphic quadratic differentials on $X$ is a $6g - 6$ dimensional vector space.

In a neighborhood of every regular point $P$ of $q$ one can introduce a local parameter $w$, in terms of which $q$ is identically equal to 1. This parameter, called the *natural parameter* of $q$ near $P$, is determined by the integral

$$w = \int_X \sqrt{q(z)}dz$$

uniquely up to a transformation $w \to \pm w + \text{const}$.

The vertical trajectories of $q$ are the arcs along which $q(z)dz^2 < 0$, and the horizontal trajectories are the arcs where $q(z)dz^2 > 0$. Hence every quadratic differential $q$ determines a pair of transverse measured foliations: the *vertical* foliation, where the leaves are the vertical trajectories, and the *horizontal* foliation, with leaves being the horizontal trajectories of $q$. One can also obtain the vertical and horizontal foliations by pulling back the vertical and horizontal foliations of $\mathbb{C}$ via a natural parameter. Every quadratic differential $q$ determines a flat metric with the length element $|q(z)|^{1/2}|dz|$. Again, the flat metric can be obtained from the natural parameter of $q$ by pulling back the Euclidean metric from $\mathbb{C}$.

Geodesic rays can be described as follows. Each direction at a point $X$ in $\mathcal{T}_g$ is associated to a quadratic differential $q$ on $X$. For $t \in \mathbb{R}$, let $q_t$ be 1–parameter family of quadratic differentials obtained from $q$ so that if $z = x + iy$ are natural coordinates for $q$ away from zeroes then $z_t = e^{-t/2}x + ie^{t/2}y$ are natural coordinates for $q_t$. Let $X_t$ be the conformal structure corresponding to $q_t$. Then $\{X_t\}$ is a geodesic.

Let $\mathcal{S}$ be the set of homotopy classes of essential simple closed curves on $M$ with the discrete topology. Here by essential we mean homotopically nontrivial and non-peripheral. The space of functionals $R_+^{\mathcal{S}}$ is given the product topology. Let $PR_+^{\mathcal{S}}$ be the corresponding projective space. Recall that the Teichmüller space $\mathcal{T}_g$ can be identified with the space of equivalence classes of hyperbolic metrics $\rho$ on $M$ of constant curvature $-1$, where $\rho_1 \sim \rho_2$ if there exists an isometry from $(M, \rho_1)$ to $(M, \rho_2)$ isotopic to the identity. Thurston showed that the hyperbolic lengths of curves can be approximated, as one goes to infinity in $\mathcal{T}_g$, by their intersection numbers with a measured foliation. Therefore one can define a compactification in terms of ratios of hyperbolic lengths. The map $\mathcal{T}_g \mapsto PR_+^{\mathcal{S}}$ defined by $\rho \to (\gamma \to \ell_\rho(\gamma))$, where $\ell_\rho(\gamma)$ is the length of the unique geodesic in the hyperbolic metric in the class of $\gamma$, is injective. It is called the Thurston embedding of Teichmüller space. The boundary of $\mathcal{T}_g$ in $PR_+^{\mathcal{S}}$ is the sphere $\mathcal{PMF}$ of projective measured foliations on $M$. The union of $\mathcal{T}_g$ and $\mathcal{PMF}$ is denoted by $\widehat{\mathcal{T}_g}$ and is called Thurston compactification. It is homeomorphic to a closed $6g - 6$ dimensional ball, where $\mathcal{PMF}$ is the $6g - 7$ dimensional sphere.

## 3  Construction and idea of the proof

Start with a square torus, ie the unit square with lower left vertex at $(0,0)$, with opposite sides identified. Cut the torus along a line segment (call it a slit) of length $0 < s < 1$ and a slope $\theta_1 > 0$. Take a second copy of the torus, cut it along a slit of same length $s$ and a slope $\theta_2 > 0$. Rotate the squares counterclockwise so that the slits are vertical. Now identify the left side of the slit on one copy with the right side on the other copy. This defines a flat structure with a parallel line field which corresponds to a quadratic differential $q_{\theta_1,\theta_2}$ with 2 zeroes of order 2 on a surface $X$ of genus 2. For more details on this construction, flat structures and quadratic differentials see Masur and Tabachnikov [9]. $X$ is partitioned into 2 sheets $S_1$ and $S_2$ separated from each other by the union of the two slits. If both $\theta_1$ and $\theta_2$ are irrational, then the vertical foliation of $q_{\theta_1,\theta_2}$ has one closed leaf, which is the union of the two slits, and all other leaves are dense in $S_1$ or $S_2$. Let $\{X_t\}$ be the Teichmüller geodesic ray from $X$ determined by the differential $q_{\theta_1,\theta_2}$.

Recall that each $x \in \mathbb{R}$ admits a continued fraction expansion of the form

$$x = a_0 + \cfrac{1}{a_1 + \cfrac{1}{a_2 + \cdots}},$$

with $a_0 \in \mathbb{Z}$, $a_i \in \mathbb{N}$, $i \geq 1$ and that $x \in \mathbb{R} - \mathbb{Q}$ iff infinitely many $a_i$'s are nonzero. We will also use the notation $x = [a_0; a_1, a_2, \ldots]$. We will call $a_i$'s the *elements* of $x$.

We will prove:

**Theorem 2** *Suppose $\theta_1 \in \mathbb{R} - \mathbb{Q}$ has bounded elements $a_{1,n} \geq 3$, $n \geq 1$, and $\theta_2 \in \mathbb{R} - \mathbb{Q}$ has unbounded elements $a_{2,n} \geq 3$, $n \geq 1$. Let $\{X_t\}$ be the Teichmüller geodesic ray constructed as above. Then $\{X_t\}$ does not converge in $\widehat{\mathcal{T}_2}$.*

Theorem 1 follows immediately from Theorem 2.

**Idea of the proof**  By the definition of Thurston compactification, a sequence $\{X_n\}$ in Teichmüller space converges to projective measured foliation $\mathcal{F}$ in $\mathcal{PMF}$ if there exists a sequence $\{r_n\}$ such that $r_n \to 0$, such that for any nontrivial simple closed curve $\alpha$ we have $r_n \cdot \ell_n(\alpha) \to i(\alpha, F)$ as $n \to \infty$. Here $\ell_n(\alpha)$ is the length of the shortest curve in the hyperbolic metric of $X_n$ which is homotopic to a curve $\alpha$, and $i(\alpha, F)$ denotes the measure of $\alpha$ with respect to a representative $F$ in $\mathcal{F}$. Hence, to prove that the geodesic does not have a limit in $\mathcal{PMF}$, it suffices to find a pair of simple closed curves $\alpha_1$ and $\alpha_2$ on $M$ and, after ruling out the possibility that $i(\alpha_j, F) = 0$ for both $j = 1$ and $j = 2$, show that the limit $\lim_{t \to \infty} \ell_t(\alpha_1)/\ell_t(\alpha_2)$ does not exist. In our proof $\alpha_j$ is the curve represented by the vector $(1, 0)$ before the rotation on $S_i$, $j = 1, 2$.

To prove the theorem, we need to be able to estimate hyperbolic lengths of simple closed curves on $X_t$ for large values of $t$. The geodesic ray $\{X_t\}$ is constructed using the flat metric. However, in general there is no easy way to get a good estimate of the hyperbolic length of a curve. Our method is to identify the shortest and the second shortest curves in the flat metric on each torus (we think of the surface $X_t$ as a surface glued out of two tori $S_1$ and $S_2$), and estimate their hyperbolic length. This information, together with the notion of intersection number, makes it possible to find good lower and upper bounds on the lengths of the curves $\alpha_1$ and $\alpha_2$.

The idea behind the choice of $\theta_1$ and $\theta_2$ is as follows. We will show (Lemma 1) that at any time the shortest curve on $S_i$ is a curve whose slope is a convergent of $\theta_i$. Moreover, we will prove that the hyperbolic length of the shortest curve depends on the elements of $\theta_i$. More precisely, the smallest hyperbolic length of the curve whose slope is $n$–th convergent of $\theta_i$ is roughly $1/a_{i,n+1}$, where $a_{i,n+1}$ is the $(n+1)$–st element of $\theta_i$ (Lemma 3). Hence assuming that elements of $\theta_2$ are unbounded and those of $\theta_1$ are bounded we easily find a sequence of times when the hyperbolic length of the shortest curve on $S_2$ goes to 0, while on the other side the length of the shortest curve stays bounded below. By the Collar Lemma, each time one intersects an extremely short curve, one has to cross a collar of width approximately $\log(1/\ell(\text{shortcurve}))$. Hence one would think that the curve $\alpha_2$ must become exceedingly long compared to the curve $\alpha_1$. However, we will see that the curve $\alpha_1$ intersects a curve of bounded length on $S_1$ significantly more than the curve $\alpha_2$ crosses the long collar in $S_2$. As a result, the ratio of the hyperbolic lengths along the sequence becomes arbitrarily large. On the other hand, after a shortest curve on $S_2$ reaches its minimal length, it has to grow. As its length becomes compatible with the length of the next short curve, the curve $\alpha_2$ grows faster and catches up with $\alpha_1$. At that time the ratio of the hyperbolic lengths is bounded above. We then conclude that the limit $\lim_{t \to \infty} \ell_t(\alpha_1)/\ell_t(\alpha_2)$ does not exist.

The main result is proven in Section 6. In Section 4 we consider the shortest curves in the flat metric. Estimates of the hyperbolic lengths of these curves are made in Section 5.

Section 7 is dedicated to the limit set of the geodesic ray. More precisely, in Section 7 we show that it a connected subset of one-dimensional simplex $L_{\theta_1,\theta_2}$. Each point in $L_{\theta_1,\theta_2}$ is a projective class of measured foliations topologically equivalent to the vertical foliation of $q_{\theta_1,\theta_2}$, with possibly different weights on $S_1$ and $S_2$.

**Acknowledgements** The author would like to thank her thesis advisor Howard Masur for his excellent guidance and support. Also the author is grateful to Juan Souto for multiple helpful suggestions, as well as to the referee for many useful comments.

## 4 Continued fractions and the flat metric on a torus

In this section we focus on one torus only, ignoring the slit and the other torus. Let $\theta = [a_0; a_1, a_2, a_3, \ldots]$ be any positive irrational number, with $a_i \geq 3$. The $k$–th convergent of $\theta$ is the reduced fraction

$$\frac{p_k}{q_k} = a_0 + \cfrac{1}{a_1 + \cfrac{1}{\cdots + \cfrac{1}{a_k}}}.$$

A few standard facts (see for example Khinchin [5]) about continued fractions are:

(1a) $\qquad p_{n+1} = a_{n+1} p_n + p_{n-1}$ and $q_{n+1} = a_{n+1} q_n + q_{n-1}$

(1b) $\qquad \dfrac{1}{q_n + q_{n+1}} \leq |p_n - \theta q_n| \leq \dfrac{1}{q_{n+1}}$

(1c) $\qquad \dfrac{p_{2n}}{q_{2n}} \nearrow \theta$ and $\dfrac{p_{2n+1}}{q_{2n+1}} \searrow \theta$

(1d) $\qquad |p_{n+1} q_n - q_{n+1} p_n| = 1$

Consider a standard lattice $\mathbb{Z}^2$ in $\mathbb{R}^2$. Let $g_t^\theta$ be a map given by a matrix

$$\frac{1}{\sqrt{1+\theta^2}} \begin{pmatrix} \theta e^{t/2} & -e^{t/2} \\ e^{-t/2} & \theta e^{-t/2} \end{pmatrix}$$

which is a rotation by an angle of $\frac{\pi}{2} - \tan^{-1}\theta$, followed by a horizontal stretch by a factor of $e^{t/2}$ and vertical contraction by $e^{t/2}$. For every $t$ we get a new lattice in $\mathbb{R}^2$. We will refer to the image of any vector $(q, p) \in \mathbb{Z}^2$ under the map $g_t^\theta$ as $(q, p)$–vector or $(q, p)$–curve at time $t$.

**Notation** To simplify our presentation we use $\approx, \Theta$ and $O$ defined as follows: for two sequences $x_n > 0, y_n > 0$, $x_n \approx y_n$ means $x_n/y_n \to 1$ as $n \to \infty$; $x_n = O(y_n)$ iff $\sup x_n/y_n < \infty$; $x_n = \Theta(y_n)$ iff $x_n = O(y_n)$ and $y_n = O(x_n)$.

**Lemma 1** *Suppose a $(q, p)$–vector is the shortest vector in Euclidean length at some time $t$. Suppose $t \geq \log\big((1 + a_0\theta)/(\theta - a_0)\big)$. Then $p/q$ is a convergent for $\theta$, ie* $p = p_n$ and $q = q_n$ for some $n$. Also at time $T_n = \log\left((p_n\theta + q_n)/|q_n\theta - p_n|\right)$, the $(q_n, p_n)$–vector is the shortest. For $t \in [T_n, T_{n+1}]$ the shortest vector is either $(q_n, p_n)$ or $(q_{n+1}, p_{n+1})$.

**Proof** By taking reciprocals if needed we can assume that $p/q < \theta$. Suppose first that $p_0/q_0 < p/q$. If $p/q$ is not a convergent of $\theta$, then there is a unique $n$ such that

(2) $$\frac{p_n}{q_n} < \frac{p}{q} < \frac{p_{n+2}}{q_{n+2}}.$$

We claim that the $(q_{n+1}, p_{n+1})$–vector is always shorter than the vector $(q, p)$. At time $t$ the image of $(q, p)$ is

$$g_t^\theta(q, p) = \frac{q + \theta p}{e^{t/2}\sqrt{1 + \theta^2}}(0, 1) + \frac{e^{t/2}(q\theta - p)}{\sqrt{1 + \theta^2}}(1, 0).$$

The Euclidean length of $g_t^\theta(q, p)$ satisfies

$$l_t((q, p))^2 = \frac{1}{1 + \theta^2}\left(\frac{(q + \theta p)^2}{e^t} + e^t(p - \theta q)^2\right).$$

Since $p_{n+1}q_n - q_{n+1}p_n = 1$ and (2) implies

(3) $$\frac{p_n}{q_n} < \frac{p}{q} < \frac{p_{n+1}}{q_{n+1}},$$

we can show that $p \geq p_n + p_{n+1}$ and $q \geq q_n + q_{n+1}$. Indeed, (3) implies

$$pq_n - p_nq \geq 1$$
and $$p_{n+1}q - pq_{n+1} \geq 1.$$

Multiplying the first inequality by $p_{n+1}$ and the second by $p_n$, and adding them we obtain

$$p_{n+1}q_np - p_nq_{n+1}p \geq p_n + p_{n+1}.$$

Therefore $p \geq p_n + p_{n+1}$. Similar argument shows that $q \geq q_n + q_{n+1}$.

Going back to the proof of the lemma, we see that

$$q + \theta p > q_{n+1} + \theta p_{n+1}.$$

Also $\quad |p - \theta q| = q|\frac{p}{q} - \theta| > q|\frac{p}{q} - \frac{p_{n+2}}{q_{n+2}}| > \frac{1}{q_{n+2}} \geq |p_{n+1} - \theta q_{n+1}|.$

The claim now follows from these inequalities and we conclude that $p/q$ is a convergent of $\theta$. Now suppose $p_0/q_0 > p/q$. Then we compute that $l_t((q, p))^2 > l_t((q_0, p_0))^2$ if $t \geq \log\left((1+a_0\theta)(\theta-a_0)\right)$. Hence taking $t$ large enough guarantees that the claim holds.

The function $l_t((q_n, p_n))^2$ reaches its minimum at $T_n = \log\left((p_n\theta + q_n)/|q_n\theta - p_n|\right)$ and
$$l^2_{T_n}((q_n, p_n)) = 2\frac{(q_n + p_n\theta)|q_n\theta - p_n|}{1+\theta^2}.$$

We have for all $k < n$
$$l_{T_n}((q_k, p_k))^2 > \frac{(q_k\theta - p_k)^2 e^{T_n}}{1+\theta^2} = \frac{(q_k\theta - p_k)^2(p_n\theta + q_n)}{(1+\theta^2)|q_n\theta - p_n|}$$
$$> 4\frac{(p_n\theta + q_n)}{q_{n+1}(1+\theta^2)} \geq 4\frac{(p_n\theta + q_n)|q_n\theta - p_n|}{1+\theta^2}.$$

In the estimates above we used (1a), (1b) and the assumption that the elements of $\theta$ satisfy $a_k > 2$. Similarly, for all $k > n$
$$l_{T_n}((q_k, p_k))^2 > \frac{(q_k + p_k\theta)^2}{(1+\theta^2)e^{T_n}} = \frac{(q_k + p_k\theta)^2|q_n\theta - p_n|}{(1+\theta^2)(p_n\theta + q_n)}$$
$$> 4\frac{(q_n + p_n\theta)^2|q_n\theta - p_n|}{(1+\theta^2)(p_n\theta + q_n)}$$
$$= 4\frac{(q_n + p_n\theta)|q_n\theta - p_n|}{1+\theta^2}.$$

We conclude that $l_{T_n}((q_n, p_n)) < \min\{l_{T_n}((q_k, p_k)) : k \neq n\}$. To verify the last claim of the lemma, we notice that the sequence $\{T_n\}$ is strictly increasing and $T_n \to \infty$. Also we see that, for each $n$, the function $l_t^2((q_n, p_n))$ is strictly convex and that the graphs of $l_t^2((q_n, p_n))$ and $l_t^2((q_k, p_k))$ intersect once (see Figure 1). Since $n < k$ implies $l_0^2((q_n, p_n)) < l_0^2((q_k, p_k))$, we see that at time $t \in [T_n, T_{n+1}]$ only $(q_n, p_n)$ or $(q_{n+1}, p_{n+1})$ can be the shortest. □

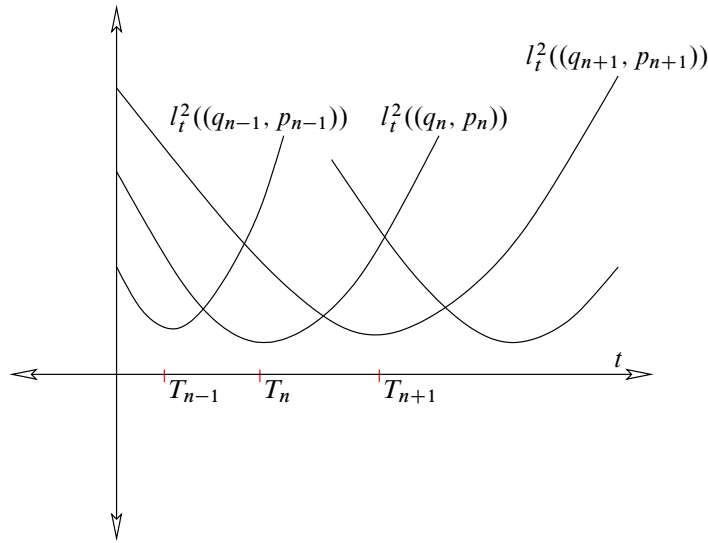

Figure 1

To find the second shortest vector is a more complicated task. There are more curves which can become second shortest than those whose slopes are convergents of $\theta$.

**Lemma 2** *Suppose that a $(q, p)$–curve is second shortest at time $t$, and $(q_n, p_n)$ is the shortest. Then either $p/q$ is a convergent of $\theta$ (namely, $(q_{n-1}, p_{n-1})$ or $(q_{n+1}, p_{n+1})$) or the following holds:*

$$q_{n-1} + q_n \leq q < q_{n+1}$$
and
$$p_{n-1} + p_n \leq p < p_{n+1}$$

**Proof** Suppose $p/q$ is not a convergent of $\theta$. Assume that $p/q < \theta$. Then there is a (unique) $k$ such that $p_{k-1}/q_{k-1} < p/q < p_{k+1}/q_{k+1}$. Using the argument from Lemma 1 we can show that the $(q_k, p_k)$–curve is shorter then $(q, p)$. Therefore $k = n$. Since the shortest and the second shortest curves intersect once, $p_n q - p q_n = 1$. Also $p_n q_{n-1} - q_n p_{n-1} = 1$. Putting these two things together we get

$$p_n(q - q_{n-1}) - q_n(p - p_{n-1}) = 0.$$

This implies that $p_n/q_n = (p - p_{n-1})/(q - q_{n-1})$. Hence $p = ap_n + p_{n-1}$ and $q = aq_n + q_{n-1}$, where $a \in \mathbb{Q}$. It is easy to see that $a < a_{n+1}$: if $a \geq a_{n+1}$ then

$$\frac{p}{q} - \frac{p_{n+1}}{q_{n+1}} = \frac{ap_n + p_{n-1}}{aq_n + q_{n-1}} - \frac{a_{n+1}p_n + p_{n-1)}}{a_{n+1}q_n + q_{n-1}} = \frac{(a_{n+1} - a)(p_{n-1}q_n - q_{n-1}p_n)}{(aq_n + q_{n-1})(a_{n+1}q_n + q_{n-1})} \geq 0$$

which is impossible. Therefore $p < p_{n+1}$ and $q < q_{n+1}$. Also $p_{n-1}/q_{n-1} < p/q < p_n/q_n$ and $p_n q_{n-1} - q_n p_{n-1} = 1$ imply $p \geq p_{n-1} + p_n$ and $q \geq q_{n-1} + q_n$. A similar argument works if $p/q > \theta$.

If $p/q$ is a convergent of $\theta$, then Lemma 1 implies that $(q, p)$ is either $(q_{n-1}, p_{n-1})$ or $(q_{n+1}, p_{n+1})$. □

## 5 The comparison of the flat and the hyperbolic metrics on the branched cover

In the previous section we considered shortest vectors on a torus. It follows from Lemma 1 that the vectors on $S_i$ with smallest length in the flat metric are those whose slopes are $n$–th convergents of $\theta_i$, which we denote by $p_{i,n}/q_{i,n}$. We claim that the same result holds for the flat surface glued from the two tori $S_1$ and $S_2$, each with a slit $(0, s)$. Indeed, since the slope of the $(p_{i,n}, q_{i,n})$–curve is a convergent, the curve can be made disjoint from the slit.

Hence at any time we know the short curves in $S_1$ and $S_2$. We want to estimate the hyperbolic length of these curves. To do so, we first prove a claim about extremal lengths of curves in $S_i$ which get short. We will denote by $\text{Ext}_t(\alpha)$ the extremal length at time $t$ of a family of curves homotopic to a curve $\alpha$. We refer the reader to Ahlfors [1] for the basic facts concerning the extremal length.

**Lemma 3** Let $\alpha$ be a $(q_{i,n}, p_{i,n})$ -curve on $S_i$, $i = 1, 2$. Then the extremal length of $\alpha$ at time $t$ satisfies
$$\frac{\text{Ext}_t(\alpha)}{l_t^2(\alpha)} \xrightarrow[n,t\to\infty]{} 1$$

**Proof** By definition of extremal length, $\text{Ext}(\alpha) \geq \inf_{\beta \sim \alpha} \left(\int_\beta \rho|dz|\right)^2 / A(\rho)$ for any metric of the form $\rho|dz|$ where $\rho \geq 0$ is Borel measurable. Assume $\alpha$ is a $(p_{1,n}, q_{1,n})$–

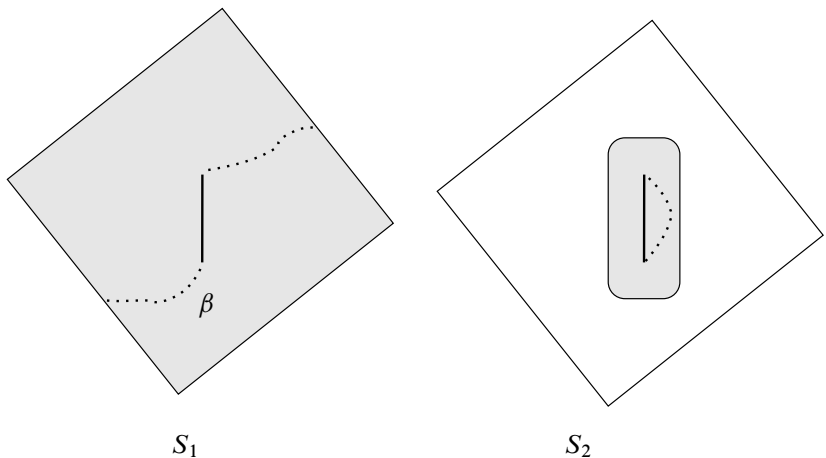

Figure 2: $\rho_t|dz|$ coincides with the flat metric in the shaded area

curve on $S_1$. Let $\rho_t|dz|$ be the metric which coincides with the flat metric of $S_1$ at time $t$, and with the flat metric of $S_2$ at time $t$ on the set of points in $S_2$ which are at most $se^{-t/2}$ (the length of the slit) away from the slit (see Figure 2). On the rest of the $S_2$ we define $\rho_t = 0$. We need to find the shortest curve with respect to this metric in the homotopy class of $\alpha$. We claim that a geodesic with respect to $\rho_t|dz|$ is contained in $S_1$, ie it does not cross the slit. The point is that any curve $\beta$ homotopic to $\alpha$ which crosses the slit will contain an arc in $S_2$ with endpoints in the slit that is homotopic relative to its endpoints into the slit. The arc is longer than the line segment connecting the endpoints. Hence there is a curve in the class of $\alpha$ which is shorter than $\beta$. Since $\alpha$ is a geodesic in the flat metric of $S_1$ at time $t$, we conclude that it is a shortest curve with respect to $\rho_t|dz|$ in its homotopy class. The length of $\alpha$ in this metric is

$$\sqrt{\frac{1}{1+\theta_1^2}\left(\frac{(q_{1,n}+\theta_1 p_{1,n})^2}{e^t}+e^t(p_{1,n}-\theta_1 q_{1,n})^2\right)}.$$

We then have

$$(4) \quad \operatorname{Ext}_t(\alpha) \geq \frac{1}{1+\theta_1^2}\frac{\frac{(q_{1,n}+\theta_1 p_{1,n})^2}{e^t}+e^t(p_{1,n}-\theta_1 q_{1,n})^2}{1+(\pi+2)s^2 e^{-t}}.$$

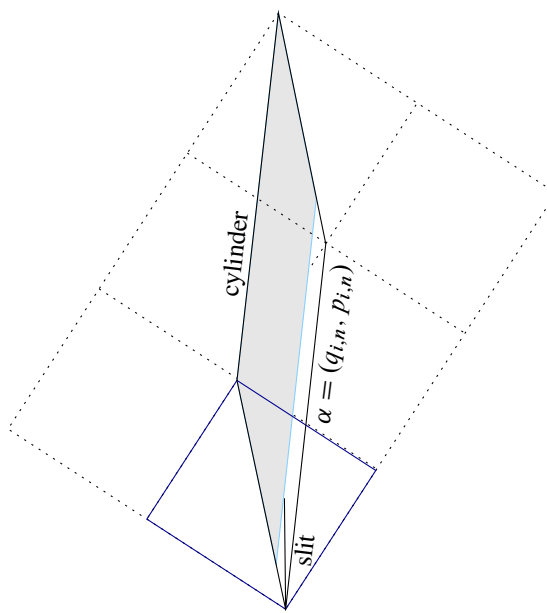

Figure 3

On the other hand, the extremal length of a simple closed curve $\alpha$ can be defined as $\inf_C(1/\operatorname{Mod}(C))$ where the infimum is taken over all cylinders $C$ with $\alpha$ a core curve, and where $\operatorname{Mod}(C)$ denotes the modulus of a cylinder $C$.

So we need to find a cylinder for a good upper bound. We consider the largest cylinder swept out by curves parallel to $\alpha$ which avoid the slit (see Figure 3). Curves parallel to $\alpha$ and crossing the slit make up a parallelogram spanned by the vectors representing the slit and $\alpha$. At time $t$ the slit is represented by a vector $u = (0, s/e^{t/2})$, and $\alpha$ is represented by

$$v = \frac{q_{1,n} + \theta p_{1,n}}{e^{t/2}\sqrt{1+\theta_1^2}}(0,1) + \frac{e^{t/2}|q_{1,n}\theta - p_{1,n}|}{\sqrt{1+\theta_1^2}}(1,0).$$

Then the area of our cylinder is

$$1 - |u \times v| = 1 - \frac{s}{e^{t/2}} \frac{e^{t/2}|q_{1,n}\theta_1 - p_{1,n}|}{\sqrt{1+\theta_1^2}} = 1 - \frac{s|q_{1,n}\theta_1 - p_{1,n}|}{\sqrt{1+\theta_1^2}}.$$

The length of the cylinder is $\|v\|$, and the height is

$$\frac{\text{area}}{\text{length}} = \frac{1 - \frac{s|q_{1,n}\theta - p_{1,n}|}{\sqrt{1+\theta_1^2}}}{\|v\|}.$$

Thus the modulus is

$$\frac{1 - \frac{s|q_{1,n}\theta - p_{1,n}|}{\sqrt{1+\theta_1^2}}}{\|v\|^2}.$$

It follows that

$$(5) \quad \operatorname{Ext}_t(\alpha) \leq \frac{1}{1+\theta_1^2}\left(\frac{(q_{1,n}+\theta_1 p_{1,n})^2}{e^t} + e^t(p_{1,n}-\theta_1 q_{1,n})^2\right)\frac{1}{1-\frac{s|q_{1,n}\theta-p_{1,n}|}{\sqrt{1+\theta_1^2}}}.$$

The claim of the lemma now follows from (5), (4) and (1b). □

**Corollary 1** Suppose $\theta_i$ is such that $a_{i,n} \to \infty$ for $i = 1$ or $i = 2$. Then at time

$$T_n = \log \frac{p_{i,n}\theta_i + q_{i,n}}{|q_{i,n}\theta_i - p_{i,n}|}$$

the hyperbolic length of a curve $\alpha = (q_{i,n}, p_{i,n})$ satisfies $\ell_{T_n}(\alpha) = \Theta(1/a_{i,n+1})$.

**Proof** By Lemma 3 we have $\operatorname{Ext}_{T_n}(\alpha) \approx l_{T_n}^2(\alpha)$. Furthermore, the assumption $a_{i,n} \to \infty$ and (1b) imply $l_{T_n}^2(\alpha) \approx (2q_{i,n})/q_{i,n+1} \approx 2/a_{i,n+1}$. On the other hand, by Proposition 1 and Corollary 3 in Maskit [7]

$$(6) \quad 2e^{-\ell(\gamma)/2} \leq \frac{\ell(\gamma)}{\operatorname{Ext}(\gamma)} \leq \pi$$

where $\text{Ext}(\gamma)$ is the extremal length and $\ell(\gamma)$ the hyperbolic length of any simple closed curve $\gamma$. In particular, if the extremal length of $\alpha$ becomes very small, then so does the hyperbolic length, and their ratio is bounded above and below. Therefore, the hyperbolic length of $\alpha$ at time $T_n$ satisfies $\ell_{T_n}(\alpha) = \Theta(1/a_{i,n+1})$. □

## 6  Proof of the main theorem

**Proof**  We begin by choosing a subsequence $\{a_{2,n_k}\}$ such that $a_{2,n_k} \to \infty$. Let for $k \in \mathbb{Z}, k \geq 1$,

$$t_{2k} = \log \frac{p_{2,n_k-1}\theta_2 + q_{2,n_k-1}}{|q_{2,n_k-1}\theta_2 - p_{2,n_k-1}|}$$

$$t_{2k+1} = \log\left((1+\theta_2^2)(q_{2,n_k})^2\right).$$

We are going to show first that the limit $\lim_{i\to\infty} \ell_{t_i}(\alpha_1)/\ell_{t_i}(\alpha_2)$ does not exist. By $\alpha_{1,1}(i)$ and $\alpha_{1,2}(i)$ (respectively $\alpha_{2,1}(i)$ and $\alpha_{2,2}(i)$) we denote curves on $S_1$ (respectively $S_2$) which are the first and second shortest in the flat metric at time $t_i$.

**i = 2k**  We will show $\sup_k \ell_{t_{2k}}(\alpha_1)/\ell_{t_{2k}}(\alpha_2) = \infty$. We have for sufficiently large $k$,

(7)  $\ell_{t_{2k}}(\alpha_2) \leq \ell_{t_{2k}}(\alpha_{2,1}(2k))i(\alpha_{2,2}(2k),\alpha_2) + \ell_{t_{2k}}(\alpha_{2,2}(2k))i(\alpha_{2,1}(2k),\alpha_2).$

Lemma 1 implies that $\alpha_{2,1}(2k)$ is a $(q_{2,n_k-1}, p_{2,n_k-1})$–curve. Thus

(8)  $i(\alpha_{2,1}(2k),\alpha_2) = p_{2,n_k-1}.$

Corollary 1 implies that

(9)  $\ell_{t_{2k}}(\alpha_{2,1}(2k)) = \Theta(1/a_{2,n_k}).$

If the second shortest curve is a $(q, p)$–curve, then by Lemma 2 we have $p \leq p_{2,n_k}$. We can assume $\alpha_{2,2}(2k)$ is the $(q_{2,n_k}, p_{2,n_k})$–curve. Therefore

(10)  $i(\alpha_{2,2}(2k),\alpha_2) = p = p_{2,n_k}.$

We need to estimate $\ell_{t_{2k}}(\alpha_{2,2}(2k))$, and for our purposes it is enough to show that

(11)  $\ell_{t_{2k}}(\alpha_{2,2}(2k)) = \Theta(\log(a_{2,n_k})).$

An easy calculation shows that the flat length of $\alpha_{2,2}(2k)$ satisfies

$$l^2_{t_{2k}}(\alpha_{2,2}(2k)) \approx a_{2,n_k}.$$

By the Lemma 3, for $t$ and $n$ large enough, the extremal length

$$\text{Ext}_{t_{2k}}(\alpha_{2,2}(2k)) \approx l^2_{t_{2k}}(\alpha_{2,2}(2k)) \approx a_{2,n_k}.$$

We will be using the well-known Collar Lemma (Theorem 4.1.1 in Buser [2]) several times in this paper. The lemma implies that any simple closed geodesic $\beta$ on a closed Riemann surface $X$ of genus $g \geq 2$ has a collar of width

(12) $$w(\beta) = \operatorname{arcsinh} \frac{1}{\sinh(1/2\ell_X(\beta))}$$

where $\ell_X(\beta)$ is the hyperbolic length of $\beta$. The collar is isometric to the cylinder

$$[-w(\beta), w(\beta)] \times \mathbb{S}^1$$

with the Riemannian metric $ds^2 = d\rho^2 + \ell_X^2(\beta) \cosh^2 \rho \, dt^2$. Hence the hyperbolic length of any curve which intersects $\beta$ nontrivially is at least $2w(\beta)$.

Let $\alpha$ be the shortest curve in the hyperbolic metric at time $t_{2k}$ such that

$$i(\alpha, \alpha_{2,2}(2k)) = 1.$$

Its length by (12) is at least

$$2 \cdot \operatorname{arcsinh}\left( \frac{1}{\frac{1}{2}\sinh(\ell_{t_{2k}}(\alpha_{2,1}(2k)))} \right) \approx \log a_{2,n_k}.$$

In fact, it is easy to see that outside the collar around $\alpha_{2,1}(2k)$ the curve $\alpha$ is bounded. Therefore $\ell(\alpha) \approx \log a_{2,n_k}$. By the estimate (6),

$$\operatorname{Ext}(\alpha) \leq \frac{1}{2}\ell(\alpha) \exp^{\ell(\alpha)/2}.$$

Minsky has shown in [10] that

$$i^2(\alpha, \alpha_{2,2}(2k)) \leq \operatorname{Ext}(\alpha) \operatorname{Ext}(\alpha_{2,2}(2k)).$$

Hence         $i(\alpha, \alpha_{2,2}(2k)) \leq \Theta(a_{2,n_k} \log^{1/2}(a_{2,n_k}))$.

Using (9), and the fact that $i(\alpha_{2,1}(2k), \alpha_{2,2}(2k)) = 1$ we see that

(13) $$\ell_{t_{2k}}(\alpha_{2,2}(2k)) = \Theta(\log a_{2,n_k})$$

and using (7), (8), (9), (10) and (13) we conclude that

(14) $$\ell_{t_{2k}}(\alpha_2) = O(\log(a_{2,n_k}) p_{2,n_k-1}).$$

Now we want to estimate $\ell_{t_{2k}}(\alpha_1)$. Recall that $\theta_1$ is such that the elements satisfy $\sup_i a_{1,i} < \infty$. Then Lemma 1 and Lemma 3 imply that $\operatorname{Ext}_{t_{2k}} \alpha_{1,1}(2k) = \Theta(1)$. Hence $\ell_{t_{2k}}(\alpha_{1,1}(2k)) = \Theta(1)$. By (12) the collar around $\alpha_{1,1}(2k)$ is of bounded width and therefore the hyperbolic length of $(\alpha_1)$ is bounded below by the number of its intersections with $\alpha_{1,1}(2k)$:

$$\ell_{t_{2k}}(\alpha_1) \geq \Theta(i(\alpha_{1,1}(2k), \alpha_1)) \tag{15}$$

If $\alpha_{1,1}(2k)$ is a $(q_{1,j_k}, p_{1,j_k})$–curve, then $i(\alpha_{1,1}(2k), \alpha_1) = p_{1,j_k} \approx q_{1,j_k}\theta_1$. Since $\alpha_{1,1}(2k)$ is the shortest curve on $S_1$ at time $t_{2k}$, it follows from Lemma 1 that the curve $(q_{1,j_k+1}, p_{1,j_k+1})$ has not reached its minimal length yet, ie

$$t_{2k} \leq \log \frac{p_{1,j_k+1}\theta_1 + q_{1,j_k+1}}{|q_{1,j_k+1}\theta_1 - p_{1,j_k+1}|}.$$

Hence for sufficiently large $k$, using (1b), we get

$$q_{2,n_k}q_{2,n_k-1}(1 + \theta_2^2) \leq \frac{p_{2,n_k-1}\theta + q_{2,n_k-1}}{|q_{2,n_k-1}\theta - p_{2,n_k-1}|}$$
$$\leq \frac{p_{1,j_k+1}\theta_1 + q_{1,j_k+1}}{|q_{1,j_k+1}\theta_1 - p_{1,j_k+1}|}$$
$$\leq (1 + \theta_1^2)q_{1,j_k+1}(q_{1,j_k+1} + q_{1,j_k+2}).$$

We then have $q_{1,j_k} \geq \Theta(q_{2,n_k-1}\sqrt{a_{2,n_k}})$. It follows that

$$\ell_{t_{2k}}(\alpha_1) \geq \Theta(q_{2,n_k-1}\sqrt{a_{2,n_k}}). \tag{16}$$

Putting together (14) and (16) we get

$$\frac{\ell_{t_{2k}}(\alpha_1)}{\ell_{t_{2k}}(\alpha_2)} \geq \frac{\Theta(q_{2,n_k-1}\sqrt{a_{2,n_k}})}{\Theta(\log(a_{2,n_k})p_{2,n_k-1})} \xrightarrow[n\to\infty]{} \infty. \tag{17}$$

**i = 2k + 1** We will demonstrate that $\sup_k \ell_{t_{2k+1}}(\alpha_1)/\ell_{t_{2k+1}}(\alpha_2)$ is bounded above. Recall that $t_{2k+1} = \log\big((1 + \theta_2^2)(q_{2,n_k})^2\big)$. Is is easy to see that $t_{2k} < t_{2k+1} < t_{2k+2}$. Hence by Lemma 1, $\alpha_{2,1}(2k)$ is either $(q_{2,n_k-1}, p_{2,n_k-1})$– or $(q_{2,n_k}, p_{2,n_k})$–curve. By Lemma 3 both $(q_{2,n_k-1}, p_{2,n_k-1})$ and $(q_{2,n_k}, p_{2,n_k})$ have at time $t_{2k+1}$ extremal (and hyperbolic) length bounded above and below. Then

$$\ell_{t_{2k+1}}(\alpha_2) \geq \Theta(p_{2,n_k}). \tag{18}$$

On the other hand $\alpha_{1,1}(2k + 1)$ and $\alpha_{1,2}(2k + 1)$ also have bounded extremal length, and therefore

$$\ell_{t_{2k+1}}(\alpha_1) = O(i(\alpha_{1,1}(2k + 1), \alpha_1)) + O(i(\alpha_{1,2}(2k + 1), \alpha_1)).$$

To estimate $i(\alpha_{1,1}(2k+1), \alpha_1)$ we argue similarly to the case when $i = 2k$. If $\alpha_{1,1}(2k+1)$ is a $(q_{1,j_k}, p_{1,j_k})$–curve, then Lemma 1 implies

$$t_{2k+1} \geq \log \frac{p_{1,j_k-1}\theta_1 + q_{1,j_k-1}}{|q_{1,j_k-1}\theta_1 - p_{j_k-1}|}.$$

Hence, using (1b), for sufficiently large $k$ we obtain

$$(1+\theta_2^2)q_{2,n_k}^2 \geq \frac{p_{1,j_k-1}\theta_1 + q_{1,j_k-1}}{|q_{1,j_k-1}\theta_1 - p_{j_k-1}|} \geq q_{1,j_k-1}(1+\theta_1^2)q_{1,j_k}.$$

Thus $i(\alpha_{1,1}(2k+1), \alpha_1) = p_{1,j_k} = \Theta(q_{1,j_k}) = O(q_{2,n_k})$. Lemma 2 implies that

$$i(\alpha_{1,2}(2k+1), \alpha_1) \leq p_{1,j_k+1} = \Theta(p_{1,j_k}).$$

We then have

(19) $$\ell_{t_{2k+1}}(\alpha_1) = O(q_{2,n_k}).$$

Putting (18) and (19) together we obtain

(20) $$\sup_n \frac{\ell_{t_{2k+1}}(\alpha_1)}{\ell_{t_{2k+1}}(\alpha_2)} < \infty.$$

To finish the proof we need to rule out the possibility that the limit is a projective measured foliation $[F]$ so that $i(\alpha_1, F) = i(\alpha_2, F) = 0$. This could be the case if the collar around $\sigma$ grew a lot faster than $\alpha_1$ and $\alpha_2$. In this case $[F] = [\sigma]$. By Theorem 3 this is impossible. $\square$

## 7 The limit set of $X_t$

In this section we say what the limit set of the geodesic ray is. Note that the only condition we put on $\theta_1$ and $\theta_2$ is that they are both irrational.

Let $\{X_t\}$ be a Teichmüller ray constructed as in Section 1, with $\theta_1, \theta_2 \in \mathbb{R} - \mathbb{Q}$. Let $F_{\theta_1,\theta_2}$ be the vertical measured foliation of the quadratic differential $q_{\theta_1,\theta_2}$. Let $\Lambda$ be the simplex of all measures on the underlying foliation of $F_{\theta_1,\theta_2}$. Also denote by $L_{\theta_1,\theta_2}$ the limit set of $\{X_t\}$ in $\mathcal{PMF}$.

**Theorem 3**  $L_{\theta_1,\theta_2}$ *is a connected subset of one dimensional simplex in* $[\Lambda]$. *If* $[F] \in L_{\theta_1,\theta_2}$ *then the weight it puts on* $\sigma$ *is 0.*

Suppose the slopes $\theta_1$ and $\theta_2$ are chosen so that one has bounded elements, and the other unbounded. It follows from Theorem 2 and Theorem 3 that $L_{\theta_1,\theta_2}$ is a one-dimensional simplex. Moreover, the proof of Theorem 2 shows that the ratio $\ell_t(\alpha_1)/\ell_t(\alpha_2)$ gets arbitrary large, and hence $L_{\theta_1,\theta_2}$ contains the projective class of the measured foliation which puts all the weight on $S_1$. Hence we have the following:

**Corollary 2** Let $\{X_t\}$ be as above. Suppose also that the elements of $\theta_1$ are bounded and those of $\theta_2$ are unbounded. Then $L_{\theta_1,\theta_2}$ is a one-dimensional simplex. Further, one of the two endpoints is such that the measure is supported on $S_1$.

**Proof of Theorem 3** Recall that a sequence $\{X_n\}$ in Teichmüller space converges to a projective measured foliation $[F]$ in $\mathcal{PMF}$ if there is $r_n \to \infty$ such that for any $\alpha \in S$ we have $r_n \ell_n(\alpha) \to i(\alpha, F)$ as $n \to \infty$. The theorem claims that the limit points are projective measured foliations which are the same as $[F_{\theta_1,\theta_2}]$, except for different ratios of weights on $S_1$ and $S_2$, which come from different accumulation points of $\{\ell_t(\alpha_1)/\ell_t(\alpha_2)\}$ as $t \to \infty$.

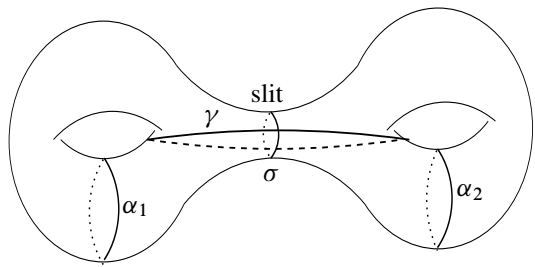

Figure 4

The fact that $L_{\theta_1,\theta_2} \in [\Lambda]$ can be demonstrated as follows. Suppose $X_{t_n} \to [F]$ as $n \to \infty$. We recall first that the space of measured foliations $\mathcal{MF}(M)$ can be identified with the space of measured laminations $\mathcal{ML}(M)$. See G Levitt's paper [6] for details. In particular, if a foliation $G \in \mathcal{MF}$ corresponds to a lamination $\mu_G \in \mathcal{ML}$ then $i(G,\cdot) = i(\mu_G,\cdot)$. Let $\mu$ and $\mu_{\theta_1,\theta_2}$ be the measured laminations corresponding to $F$ and $F_{\theta_1,\theta_2}$. Then

$$r_n \ell_{t_n}(\mu_{\theta_1,\theta_2}) \to i(\mu, \mu_{\theta_1,\theta_2}) = i(F, F_{\theta_1,\theta_2}).$$

Since $\ell_{t_n}(\mu_{\theta_1,\theta_2}) \to 0$, we conclude that $i(F, F_{\theta_1,\theta_2}) = 0$. Since $F_{\theta_1,\theta_2}$ intersects every simple closed curve except for $\sigma$, and $F$ is not minimal, $F$ has one closed leaf which is in the homotopy class of $\sigma$. It is also easy to see that any other leaf of $F$ is dense in $S_1$ or $S_2$ and has a slope of $\theta_1$ or $\theta_2$. Therefore $F \in \Lambda$.

To prove that $L_{\theta_1,\theta_2}$ is contained in the one dimensional simplex in $[\Lambda]$ we consider a curve $\gamma$ which crosses the dividing curve $\sigma$ (Figure 4) and estimate its hyperbolic length along the geodesic ray. It suffices to demonstrate that

$$(21) \qquad \frac{\ell_t(\gamma)}{\ell_t(\alpha_1) + \ell_t(\alpha_2)} \to 1.$$

At any time $t$ the image of the curve $\gamma$ can be thought of as a curve made of an arc parallel to the image of $\alpha_1$ followed by an arc crossing the slit $\sigma$ and perhaps winding around $\sigma$ for a while, then followed by an arc parallel to the image of $\alpha_2$ and then by another arc crossing $\sigma$. Hence the length of $\gamma$ in the hyperbolic metric at time t satisfies

$$\ell_t(\alpha_1) + \ell_t(\alpha_2) - C \leq \ell_t(\gamma) \leq \ell_t(\alpha_1) + \ell_t(\alpha_2) + 2c_t + C$$

where $c_t$ is the hyperbolic length of an arc connecting the geodesic (in the flat metric) representatives of $\alpha_1$ and $\alpha_2$ at time $t$ and parallel to $\gamma$, and $C$ is some positive constant independent of $t$. Therefore

$$(22) \qquad 1 - \frac{C}{\ell_t(\alpha_1) + \ell_t(\alpha_2)} \leq \frac{\ell_t(\gamma)}{\ell_t(\alpha_1) + \ell_t(\alpha_2)} \leq 1 + \frac{2c_t + C}{\ell_t(\alpha_1) + \ell_t(\alpha_2)}.$$

It suffices to show that

$$(23) \qquad \frac{c_t}{\ell_t(\alpha_1) + \ell_t(\alpha_2)} \to 0.$$

Consider the cylinder $A \in X_t$ with $\sigma$ the core curve which boundary components are Euclidean circles $C_1$ and $C_2$ of radii $R_1$ and $R_2$ (see Figure 5). We want to have an upper bound for its modulus. Since $\text{Mod}(A) = 1/\text{Ext}(\Gamma)$, where $\Gamma$ is the family of curves homotopic to a core curve of $A$, we first find a good lower bound for $\text{Ext}(\Gamma)$.

Let $z_j$ be the midpoint of the slit on $S_j$ for $j = 1, 2$. Let $A_j \subset S_j$ be the annulus centered at $z_j$, with inner radius $r_j = \frac{1}{2} s e^{-t/2}$ and outer radius $R_j$. Define a metric $\rho^t(z)$ on $A$ to be

$$\rho^t(z) = \left\{ \begin{array}{l} \frac{1}{2\pi r}, \text{ if } z \in A_j,\ j = 1, 2 \\ \frac{1}{2\pi r_j}, \text{ if } z \in A \setminus (\cup A_j) \end{array} \right\}.$$

Then it is easy to see that

$$A(\rho^t) = \frac{1}{2\pi} \log \frac{R_1 R_2}{r_1 r_2} + \frac{1}{2\pi} = \frac{1}{2\pi}(\log \frac{4 R_1 R_2}{s^2} + t + 1)$$

and $\qquad \inf_{\beta \leftrightarrow \alpha} \int_\beta \rho |dz| = \frac{1}{2\pi r_1} l_t(\sigma) = \frac{1}{2\pi s e^{-t/2}} 2 s e^{-t/2} = \frac{1}{\pi}.$

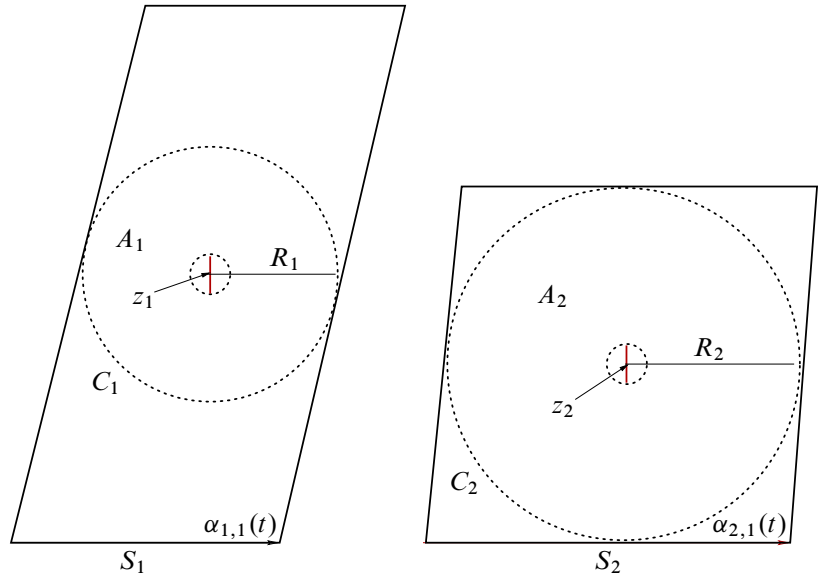

Figure 5

Hence
$$\operatorname{Ext}(\Gamma) \geq \frac{2}{\pi(\log \frac{4R_1 R_2}{s^2} + t + 1)}$$

and noticing that $R_j \leq l_t(\alpha_{j,1}(t))$ we now have

(24)  $$\operatorname{Mod}(A) \leq \frac{(\log \frac{4R_1 R_2}{s^2} + t + 1)}{2} = \Theta(t).$$

Now the cylinder is conformally equivalent to the annulus
$$\mathcal{A} = \{z \in \mathbb{C} \mid e^{-2\pi \operatorname{Mod}(A)} < |z| < 1\}$$

with the hyperbolic metric [8, Lemma 3, p 188]
$$\rho_t(z) = \frac{-\pi |dz|}{|z|(2\operatorname{Mod}(A)) \sin(\log |z|/(2\operatorname{Mod}(A)))}.$$

Further, applying Lemma 4 in [8] we see that $c_t$ is approximately equal to the length of the radius of $\mathcal{A}$. More precisely, given $\epsilon > 0$ there is a $\delta > 0$ such that

(25)  $$\left| \frac{\rho_{X_t}(z)}{\rho_t(z)} - 1 \right| < \epsilon, \, z \in \mathcal{A}_\delta = \{e^{-2\pi m}/\delta < |z| < \delta\}$$

where $\rho_{X_t}$ is the hyperbolic metric of $X_t$. A simple calculation shows that a radius of $\mathcal{A}_\delta$ in the metric $\rho_t(z)$ is

$$\left( \log \frac{1 + \cos \frac{\log 1/\delta}{2\,\mathrm{Mod}(A)}}{1 - \cos \frac{\log 1/\delta}{2\,\mathrm{Mod}(A)}} \right)(1 + \epsilon) = \Theta(\log(\mathrm{Mod}(A))).$$

Hence $$c_t = \Theta(\log(\mathrm{Mod}(A))) = O(\log t).$$

On the other hand, if $\alpha_{j,1}(t)$ is a $(q_{j,n_j}, p_{j,n_j})$ curve for some $n_j$ and $j = 1, 2$ then

$$e^{t/2} = O\left( \frac{1}{|q_{j,n_j} - \theta_j \, p_{j,n_j}|} \right) = O(q_{j,n_j} + q_{j,n_j+1})$$

and hence $\quad t = O(\log(q_{j,n_j} + q_{j,n_j+1})) = O(\log a_{j,n_j+1} q_{j,n_j})$.

Therefore we have $\quad c_t = O(\log \log a_{j,n_j+1} q_{j,n_j}))$.

Since $\ell_t(\alpha_j) \geq \Theta(p_{j,n_j} \log(a_{j,n_j+1}))$, it is clear that (23) holds.

The function $\ell_t(\alpha_1)/\ell_t(\alpha_2)$ is continuous, and therefore the limit set is an interval or a point. This concludes the proof. $\square$

# References

[1] **L V Ahlfors**, *Lectures on quasiconformal mappings*, second edition, University Lecture Series 38, American Mathematical Society, Providence, RI (2006)   MR2241787With supplemental chapters by C. J. Earle, I. Kra, M. Shishikura and J. H. Hubbard

[2] **P Buser**, *Geometry and spectra of compact Riemann surfaces*, Progress in Mathematics 106, Birkhäuser, Boston (1992)   MR1183224

[3] **A Fathi**, **F Laudenbach**, **V Poenaru**, *Travaux de Thurston sur les surfaces*, Astérisque 66, Société Mathématique de France, Paris (1979)   MR568308 Séminaire Orsay, with an English summary

[4] **Y Imayoshi**, **M Taniguchi**, *An introduction to Teichmüller spaces*, Springer, Tokyo (1992)   MR1215481Translated and revised from the Japanese by the authors

[5] **A Y Khinchin**, *Continued fractions*, russian edition, Dover Publications, Mineola, NY (1997)   MR1451873With a preface by B. V. Gnedenko, Reprint of the 1964 translation

[6] **G Levitt**, *Foliations and laminations on hyperbolic surfaces*, Topology 22 (1983) 119–135   MR683752

[7] **B Maskit**, *Comparison of hyperbolic and extremal lengths*, Ann. Acad. Sci. Fenn. Ser. A I Math. 10 (1985) 381–386   MR802500

- [8] **H Masur**, *Two boundaries of Teichmüller space*, Duke Math. J. 49 (1982) 183–190 MR650376
- [9] **H Masur**, **S Tabachnikov**, *Rational billiards and flat structures*, from: "Handbook of dynamical systems, Vol. 1A", (B Hasselblatt, A Katok, editors), North-Holland, Amsterdam (2002) 1015–1089   MR1928530
- [10] **Y N Minsky**, *Teichmüller geodesics and ends of hyperbolic* 3*-manifolds*, Topology 32 (1993) 625–647   MR1231968
- [11] **K Strebel**, *Quadratic differentials*, Ergebnisse der Mathematik und ihrer Grenzgebiete (3) 5, Springer, Berlin (1984)   MR743423

*Department of Mathematics, University of Michigan*
*530 Church Street, Ann Arbor MI 48109, USA*

alenzhen@umich.edu